\documentclass[12pt, dvips, twoside]{amsart}
\usepackage{amsmath,amssymb,amscd,amsxtra,latexsym,bbm,graphicx,epsfig,
oldgerm,bm,euscript,eufrak,psfrag,amsthm}
\usepackage[active]{srcltx}
\usepackage{a4wide}

\theoremstyle{plain}
\newtheorem{theorem}{Theorem}[section]
\newtheorem*{theorem*}{Theorem}
\newtheorem{lemma}[theorem]{Lemma}

\newtheorem*{remark*}{Remark}
\newtheorem*{remarks*}{Remarks}

\newtheorem{remarks}[theorem]{Remarks}

\newtheorem*{example*}{Example}
\newtheorem*{examples*}{Examples}

\newtheorem*{definition*}{Definition}

\newcommand*{\lou}[1]{{\raisebox{-.2ex}{\mathsurround=0pt\makebox{$\scriptstyle\mskip-.8mu#1$}}}}

\newcommand{\proofend}{\hspace*{\fill} $\Box$\\}

\def\s{\smallskip}
\def\m{\medskip}

\def\Diffc0{\operatorname{Diff^c_0}}

\def\Sympc0{\operatorname{Symp^c_0}}

\def\HP{\operatorname{HP}}

\def\e{\operatorname{e}}
\def\i{\operatorname{i}}
\def\aff{\operatorname{aff}}

\def\ga{\alpha}
\def\gb{\beta}
\def\gg{\gamma}

\def\eps{\varepsilon}

\def\gf{\varphi}

\def\gl{\lambda}
\def\go{\omega}
\def\gs{\sigma}

\def\cb{{\mathcal B}}

\def\cf{{\mathcal F}}

\def\ch{{\mathcal H}}
\def\cj{{\mathcal J}}

\def\cl{{\mathcal L}}

\def\cs{{\mathcal S}}
\def\ct{{\mathcal T}}

\def\T{\widehat{\mbox{$\mskip-.2mu T\mskip1mu$}}}

\def\CC{\mathbbm{C}}

\def\NN{\mathbbm{N}}

\def\RR{\mathbbm{R}}

\def\ZZ{\mathbbm{Z}}

\def\CP{\operatorname{\mathbbm{C}P}}

\def\pp{\partial}

\def\iiota{\dot\iota}

\def\ni{\noindent}
\def\b{\bigskip}
\def\m{\medskip}

\def\.{\mskip1mu}
\def\?{\mskip-1mu}

\def\id{\operatorname{\mbox{\tt id}}}

\def\e{\operatorname{e}}
\def\i{\operatorname{i}}

\begin{document}

\title{Notes on monotone Lagrangian twist tori}

\author{Yuri Chekanov}
\address{(Yu.~Chekanov)
Moscow Center for Continuous Mathematical Education,
B.~Vlasievsky per.~11, Moscow 121002, Russia}
\email{chekanov@mccme.ru}
\author{Felix Schlenk} \thanks{FS partially supported by SNF grant 200021-125352/1.}
\address{(F.~Schlenk)
Institut de Math\'ematiques,
Universit\'e de Neuch\^atel,
Rue \'Emile Argand~11,
CP~158,
2009 Neuch\^atel,
Switzerland}
\email{schlenk@unine.ch}

\date{\today}

\begin{abstract}
We construct monotone Lagrangian tori 
in the standard symplectic vector space, 
in the complex projective space and in products of spheres.
We explain how to classify these Lagrangian tori 
up to symplectomorphism and Hamiltonian isotopy,
and how to show that they are not displaceable by Hamiltonian isotopies.
\end{abstract}

\maketitle

\markboth{{\rm }}{{Twist tori}}

\section{Introduction}

\ni
These are notes for the second author's talk at MSRI, March~2010,
on certain exotic monotone Lagrangian tori (called twist tori).
More results and proofs are given in~\cite{CS1, CS2}.

\m
A {\bf Lagrangian torus} in a $2n$-dimensional symplectic manifold $(M,\go)$
is a submanifold~$L$ diffeomorphic to an $n$-dimensional torus such that
$\go$ vanishes on the tangent bundle $TL$.
Given such an $L$, consider the area homomorphism $\gs \colon \pi_2(M,L) \to \RR$
defined by $\gs \bigl( [D] \bigr) = \int_D \go$ and the Maslov homomorphism
$\mu \colon \pi_2(M,L) \to \ZZ$ defined as in~\cite{Ar}.
Then $L$~is {\bf monotone} if $\gs = C \mu$ for some constant $C>0$.

\s
There are by now many strong tools to study monotone Lagrangian submanifolds,
such as Floer homology \cite{Fl,FOOO1,FOOO2}, pearl homology \cite{BC1,BC2},
and symplectic quasi-states~\cite{EP}.
Except for fibers of toric symplectic manifolds, however, only a few
examples of monotone Lagrangian tori are known.
We provide many such examples.
They can be used to test and refine the existing tools.

\section{Construction}

\ni
Fix $k \in \NN$ and $\eps >0$.
Denote by $D^2(a)$ the open disc of area $a$ in $\RR^2$
centred at the origin, and by $\cs(k)$ the open sector in $\RR^2$,
$$
\cs (k) \,=\, \left\{ r \e^{\i \gf} \mid 0<\gf< \frac{2\pi}{k+1}\right\}.
$$
Let $\gg$ be a smooth embedded curve in $\RR^2$
such that
\begin{itemize}
\item[$\bullet$]\;
$\gg$ encloses a domain of area $1$;
\item[$\bullet$]\;
$\gg$ lies in the sector $\cs(k) \cap D^2(k+1+\eps)$.
\end{itemize}
\begin{figure}[ht]
 \begin{center}
 \psfrag{k}{$\frac{2\pi}{k+1}$}
  \leavevmode\epsfbox{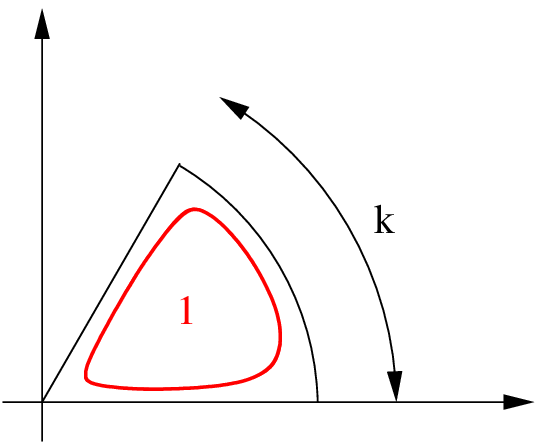}
 \end{center}
 \caption{The curve $\gg$ for $k=5$.}
\end{figure}
%
%

\begin{figure}[ht]
 \begin{center}
  \leavevmode\epsfbox{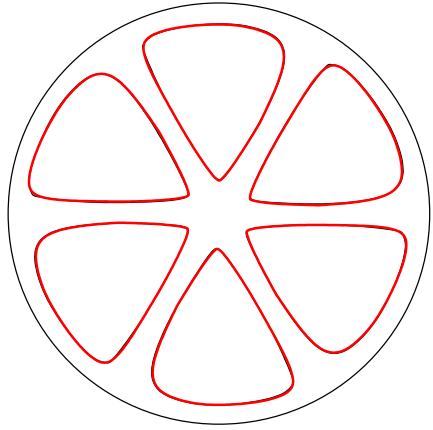}
 \end{center}
 \caption{$\Theta^5 \cap \Delta$.}
 \label{fig.rose}
\end{figure}

\ni
Define the ``twist torus'' $\Theta^k$ in $\RR^{2(k+1)}$ as
\begin{eqnarray*}
\Theta^k &=&
\left\{ \frac 1{\sqrt{k+1}} \Bigl( \e^{\i \ga_1} \gg(t), \,\e^{\i \ga_2} \gg(t), \,\dots,
           \,\e^{\i \ga_{k+1}} \gg(t) \Bigr)  \;\Big|\;  \sum_{j=1}^{k+1} \ga_j = 0 \right\}
\end{eqnarray*}
where $\gg(t)$ is a parametrization of $\gg$, and $\ga_j \in \RR$.
The torus
\begin{eqnarray*}
\Theta^1 &=&
\left\{ \frac 1{\sqrt 2} \Bigl( \e^{\i \ga} \gg(t), \e^{-\i \ga} \gg(t) \Bigr)  \right\}
\end{eqnarray*}
in $\RR^4$ was first described in \cite{Ch-tori} and \cite{ElP}.
Note that the curve $\Gamma (t) := \frac 1{\sqrt{k+1}} \bigl( \gg(t), \dots, \gg(t) \bigr)$
lies in the diagonal plane
$\Delta = \left\{ (z,\dots,z) \right\} \subset \CC^{k+1} \cong \RR^{2(k+1)}$ and that
$\Theta^k$ is obtained by restricting the action
$(z_1,\dots,z_{k+1}) \mapsto \left( \e^{\i \ga_1} z_1, \,\dots, \, \e^{\i \ga_{k+1}} z_{k+1} \right)$
of $\mathrm{T}^{k+1}$ on $\CC^{k+1}$ to the $k$-dimensional subtorus
$$
\mathrm{T}_0^k \,=\, \biggl\{ \bigl( \e^{\i \ga_1} , \, \dots, \, \e^{\i \ga_{k+1}} \bigr)
                                      \;\Big|\;  \sum_{j=1}^{k+1} \ga_j = 0 \bigg\}
$$
and to the curve $\Gamma$.
%
%
The torus $\Theta^k$ is embedded because $\gg$ lies in the sector
$\cs(k)$.
Indeed, the intersection $\Theta^k \cap \Delta$ looks as in Figure~\ref{fig.rose}.

%
%

\ni

The torus $\Theta^k$ is Lagrangian because the orbits of the group~$\mathrm{T}_0^k$
are isotropic and skew orthogonal to~$\Delta$.
It is monotone in $\RR^{2(k+1)}$ since its Maslov class and
symplectic area class vanish on the orbits of~$\mathrm{T}_0^k$.
Since $\gamma$ lies in $D^2(k+1+\eps)$,
the torus $\Theta^k$ belongs to the polydisc

\begin{equation} \label{e:poly}
D^{2(k+1)}(1+\eps)\, \:=\, D^2(1+\eps) \times \dots \times D^2(1+\eps).
\end{equation}

\m \ni
{\bf Let's twist again\;\!!}
This twisting construction can be iterated as follows.
For $k \in \NN$ and $z \in \cs (k)$, define the isotropic $k$-torus in $\CC^{k+1}$
by
$$
\Theta^k (z) \,=\,
\left\{ \frac 1{\sqrt{k+1}} \Bigl( \e^{\i \ga_1} z,\,\dots,\, \e^{\i \ga_{k+1}} z \Bigr)
\;\Big|\;  \sum_{j=1}^{k+1} \ga_j = 0 \right\} .
$$
Consider now a subset
$A \subset D^{2(k_1+1)}(1+\eps)$,
and take $k_2 \in \NN$ and $\ell_2 \in \{ 1, \dots, k_1+1 \}$.
The number $\ell_2$ designates the coordinate plane of $\CC^{k_1}$
on which we perform the twist.
For notational convenience we assume $\ell_2=1$.
Choose a symplectic embedding $\psi \colon D(1+\eps) \to \cs (k_2) \cap D(k_2+1+\eps)$.
Define
\begin{eqnarray*}
\Theta^{k_2}_{\ell_2} (A) &=&
\left\{ \Theta^{k_2}(z_1) \times (z_2, \dots, z_{k_1+1}) \mid (z_1, z_2,\dots, z_{k_1+1})
\in \left( \psi \times \id_{k_1} \right)(A) \right\} .
\end{eqnarray*}
Note that for a curve $\gg$ as above, $\Theta^k_1(\gg) = \Theta^k$.
If $A$ is a submanifold in $D^{2(k_1+1)}(1+\eps)$, then $\Theta^{k_2}_{\ell_2} (A)$ is a submanifold
of $D^{2(k_1+k_2+1)}(1+\eps)$ diffeomorphic to $T^{k_2} \times A$,
and if $A$ is a (monotone) Lagrangian submanifold then so is $\Theta^{k_2}_{\ell_2} (A)$.
In particular,
$\Theta^{k_2}_{\ell_2} \bigl( \Theta^{k_1} \bigr) =: \Theta^{k_2}_{\ell_2} \Theta^{k_1}$
is a (monotone) Lagrangian torus in $D^{2(k_1+k_2+1)}(1+\eps)$.
Iterating this construction, we obtain for natural numbers $k_1, \dots, k_m$ and
$1 = \ell_1, \ell_2, \dots, \ell_m$ with $1 \le \ell_j \le k_1 + \dots + k_{j-1}+1$
a monotone Lagrangian torus
\begin{equation} \label{e:subset}
\Theta_{\bf \ell}^{\bf k} \,:=\,\Theta_{\ell_m}^{k_m}
\cdots \Theta_{\ell_2}^{k_2} \Theta^{k_1}_{\ell_1} \,\subset\, D^{k_1+\dots+k_m+1}(1+\eps) .
\end{equation}
We call a Lagrangian torus of the form $\Theta_{\ell}^{\bf k}$ a {\bf primitive twist torus}.
By a {\bf twist torus} we understand a product of primitive twist tori.

\b \ni
{\bf Graphical representation.}
We represent the tori $\Theta^{\bf k}_{\ell}$ and their products by planar rooted forests.
To a primitive twist torus $\Theta^{\bf k}_{\ell} = \Theta_{\ell_m}^{k_m} \cdots \Theta_{\ell_2}^{k_2} \Theta^{k_1}_{\ell_1}$
we associate a rooted tree $\ct \bigl( \Theta_{\ell}^{\bf k} \bigr)$ recursively as follows.
To the circle in $\RR^2$ we associate a point (the root), and to $\Theta^k$
we associate the bush $\cb (k+1)$ with $k+1$ leaves connected to the root.
With $\Theta^{k_m}_{\ell_m} \left(\Theta^{\bf k}_{\ell} \right)$
we associate the tree obtained by gluing the root of the bush $\cb(k_m+1)$ to the
$\ell_m$'th leaf (counted from the left) of the tree $\ct \left( \Theta_\ell^{\bf k} \right)$.
The planar rooted forest $\cf$ associated to a product of primitive twist tori
is the disjoint union of the rooted trees associated to the factors.
See the figure below for examples.

A rooted tree is called {\bf ample}
if it is a point or if the valency at the root is at least two and the valency at each vertex that is neither a root nor a leaf is at least three.
The set of rooted trees associated to primitive twist tori of dimension~$n$
is exactly the set of ample rooted trees with~$n$ leaves.

\begin{figure}[ht]
 \begin{center}
  \psfrag{1}{$S^1 \approx \gamma \subset \RR^2 : $}
  \psfrag{2}{$\Theta^1 \subset \RR^4 : $}
  \psfrag{3}{$\Theta^4 \subset \RR^{10} : $}
  \psfrag{4}{$\Theta_1^1 \Theta^1 : $}
  \psfrag{5}{$\Theta_3^2 \Theta_1^1 \Theta^2 : $}
  \psfrag{6}{$\Theta_2 \times \Theta_1 \Theta^1 : $}
  \psfrag{7}{$T^4_{\operatorname{Cliff}} : $}
  \leavevmode\epsfbox{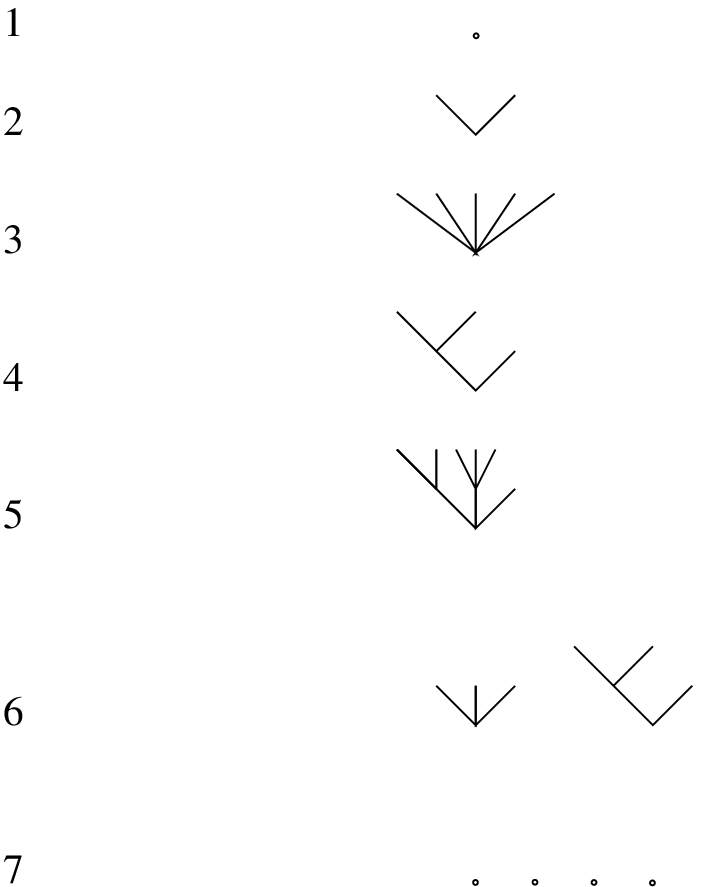}
 \end{center}
 \label{fig.ex}
\end{figure}
%
%

\section{Classification}

Two Lagrangian submanifolds $L_1, L_2$ of a symplectic manifold $(M,\go)$
are {\bf Hamiltonian isotopic} if there exists a Hamiltonian diffeomorphism
of $(M,\go)$ that maps $L_1$ to $L_2$.
Moreover,
$L_1, L_2$ are {\bf symplectomorphic}
if there exists a symplectomorphism of $(M,\go)$ mapping $L_1$ to $L_2$.
In all our results, the Hamiltonian and the symplectic classification agree.
We shall thus restrict ourselves to classification up to symplectomorphism.

In $\RR^{2n}$ twist tori form only $n$ equivalence classes up to symplectomorphism
(and scaling),
namely those formed by the tori already found in~\cite{Ch-tori}.
If we map our twist tori to closed manifolds by a Darboux chart,
there are often many more equivalence classes, however.
To fix the ideas, we look at twist tori in products of spheres
$\times_n S^2 := S^2 \times \dots \times S^2$.
Similar results can be obtained for twist tori in complex projective spaces, their products,
and their monotone blow-ups.

Let $S^2$ be the $2$-sphere endowed with an area form of total area~$2$.
Let $\infty$ and $0$ be the north and south pole of $S^2$.
Choose a symplectomorphism $\psi \colon D^2(2) \to S^2 \setminus \infty$
such that $\psi(0)=0$.
In our construction of twist tori we can keep $\eps >0$ as small as we like.
In view of~\eqref{e:subset}, each twist torus in $\RR^{2n}$ then lies in $D^{2n}(2)$.
Under the product embedding
$$
\psi \times \dots \times \psi \colon D^{2n}(2) \to \times_n S^2
$$
twist tori are mapped to Lagrangian tori in $\times_n S^2$,
that are again called twist tori.
By putting the polydisc $D^{2n}(1+\eps)$ into the open ball $B^{2n}(n+1)$ of capacity $\pi r^2 = n+1$
and by symplectically compactifying this ball into the standard complex projective space,
we also construct twist tori in~$\CP^n$.
The size $n+1$ of the ball is chosen such that these tori in $\CP^n$ are monotone.
One can also interpolate between these two cases, $\times_n S^2$ and $\CP^n$,
by constructing twist tori in products of complex projective spaces.

By an isomorphism of rooted forests we mean a homeomorphism that maps roots to roots.

\m \ni
{\bf Theorem 1.}
{\it
Two twist tori in $\times_n S^2$ are symplectomorphic if and only if
their rooted forests are isomorphic.

For $n\le7$, two twist tori in $\CP^n$
are symplectomorphic if and only if
their rooted forests are isomorphic.
}

\m
For $n\ge8$, there are twist tori in $\CP^n$ that correspond to non-isomorphic
trees but cannot be distinguished by our methods.

Ample rooted trees (and forests) with $n$ leaves can be enumerated; their number grows like
$c^n$ where $c \approx 3.692$.
\footnote{
The sequence $a_n$ of the number of rooted trees with $n$ leaves
is the sequence A000669, ``Number of series-reduced planted trees with n leaves''
of the On-Line Encyclopedia of Integer Sequences.}
%

\m
There are two proofs of this theorem.
One of them is by enumerating $J$-holomorphic discs.
The other proof is elementary in the
sense that it does not use Floer homology or $J$-holomorphic discs.
Its main tool is the {\bf displacement energy} of a set.
Denote by $\Phi_{\?\?H}$ the time 1 map of the Hamiltonian flow generated by a smooth function
$H \colon [0,1] \times M \to \RR$.
Following \cite{Hofer90}, define the norm of~$H$ by
\[
\| H \| \,=\, \int_0^1 \left( \max_{x \in M} H(t,x) - \min_{x \in M} H(t,x) \right) dt ,
\]
and the displacement energy of a compact subset $A \subset M$ by
\[
e \left( A \right) \,=\, \inf_{H \in \ch}
\Bigl\{ \| H \|
\mid \Phi_H (A) \cap A = \varnothing \Bigr\} ,
\]
assuming $\inf(\varnothing)=\infty$.
By Theorem~2 below, the displacement energy of a twist torus in $\times_n S^2$
(or $\CC P^n$) is infinite.
At first sight this invariant therefore cannot be used to distinguish twist tori.
We look, however, at {\bf nearby} tori, following~\cite{Ch-tori}.

\m \ni
{\bf Versal deformations}.
Let $(M,\go)$ be a symplectic manifold, and
denote by $\cl$ the space of closed embedded Lagrangian submanifolds
in $(M,\go)$ endowed with the $C^\infty$-topology.
The displacement energy $e$ is a function  on $\cl$ with values in $[0,+\infty]$.
For each $L \in \cl$, it gives rise to a function germ
$S^e_L \colon H^1(L;\RR) \to [0,+\infty]$ at the point $0\in H^1(L;\RR)$,
which we call the {\bf displacement energy germ}.
By Weinstein's Lagrangian Neighbourhood Theorem,
there is a symplectomorphism from
a neighbourhood of the zero section of~$T^*\?\?L$
to a neighbourhood of $L$ in $M$
such that the zero section is mapped to~$L$.
Given a sufficiently small $\xi\in H^1(L;\RR)$,
define  $S^e_L(\xi)=e(L_\xi)$, where $L_\xi$ is the
image of a closed 1-form  on~$L$ representing the class~$\xi$.
Displacement energy germs are symplectically invariant
in the following sense: for each symplectomorphism $\psi$
we have
\begin{equation}  \label{e:Sfl}
S^e_{\psi(L)} \,=\,  S_L^e \circ \left( \psi {\mid}_{L}\right)^*.
\end{equation}

\m
For the proof of Theorem~1
we pick bases in $H^1(L;\ZZ)$ and $H^1(L';\ZZ)$ to identify these groups with $\ZZ^n$,
and show that for twist tori $L,L'$ with non-isomorphic forests the functions
$S^e_{L}$ and $S^e_{L'}$ are not related by a linear isomorphism of $\ZZ^n$,
whence $L,L'$ are not symplectomorphic by~\eqref{e:Sfl}.

\b
In order to compute the functions $S^e_{L}$
we need yet another invariant of a closed Lagrangian submanifold~$L$
in a closed symplectic manifold~$(M,\go)$,
namely the {\bf Gromov width of~$L$} as defined in~\cite{Gr}.
Denote by $D$ the closed unit disc in the complex plane~$\CC$,
and by $\cj=\cj(M,\go)$ the set of tame almost complex structures $J$ on~$M$.
Given $J \in \cj$, define $\gs (L,J)$ to be the minimal symplectic area $\int_D u^*\go$
of a non-constant $J$-holomorphic map
$u \colon (D, \pp D) \to (M,L)$ if such maps exist,
and set $\gs (L,J)=\infty$ otherwise.
Define
\[
\gs\left (L\right) \,=\, \sup_{J \in \cj} \.\gs (L,J) .
\]
It was proved in \cite{Ch-energy} that
\begin{equation}   \label{e:energy}
   \gs \left(L\right) \,\le\, e \left(L\right).
\end{equation}

\m
We illustrate the computation of the function $S_L^e$ for the
Clifford torus $T^2$ and the twist torus $\Theta$ in $M = S^2 \times S^2$.

\m
\ni
{\bf Computation of $S_{T^2}^e$.}
Recall that $M_{\aff} = D^2(2) \times D^2(2)$ is the affine part of $M=S^2 \times S^2$.
Consider the standard 2-torus action on $M_{\aff}$
$$
(u,v) \,\mapsto\, \left( \e^{2\pi \i \ga}u, \e^{2\pi \i \gb}v \right) .
$$
Its moment map $\mu \colon M_{\aff} \to \RR^2$ is given by
$\mu (u,v) = \pi \left( |u|^2,|v|^2  \right)$.
The moment polyhedron $\Box = \mu (M_{\aff})$ is the square
$[0,2[ \times [0,2[$,
and the Clifford torus $T^2$ is, by definition, the preimage $\mu^{-1}(b)$ of the barycentre~$b=(1,1)$ of $\Box$.
A neighbourhood of $0 \in H^1(T^2;\RR) \cong H_1(T^2;\RR)$ can be identified with a neighbourhood of~$b$ in $\Box$.

\begin{figure}[ht]
 \begin{center}
 \psfrag{1}{$1$}
 \psfrag{2}{$2$}
  \leavevmode\epsfbox{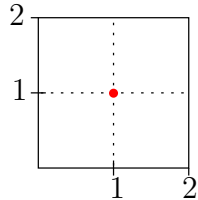}
 \end{center}
 \caption{The image $\mu(T^2) \in \Box$ of the Clifford torus.}
 \label{fig.cliff}
\end{figure}
%
%

\ni
Fix a point $p$ on the equator $E$ of $S^2$.
Take as a basis of $H_1(T^2;\ZZ)$ the classes $[E \times p]$, $[p \times E]$
of the equators,
and use this basis to identify $H_1(T^2;\RR)$ with $\RR^2$.
Then $(0,0)~\in~\RR^2$ corresponds to $T^2$, and $(x,y)$ close to $(0,0)$ corresponds to the nearby
product torus $T^2_{x,y} = \mu^{-1}(1+x,1+y)$.
For $\eps >0$, a disc in $S^2$ of area $a<1$ can be displaced from itself by a Hamiltonian isotopy
of energy smaller than $a+\eps$.
Therefore, $e(T_{x,y}) \le \min \left\{ 1-|x|, 1-|y| \right\}$ whenever $(x,y) \neq (0,0)$.
Moreover, for the standard complex structure $J_0 = \i \oplus \i$ on $S^2 \times S^2$
we have
$$
\gs \left( T_{x,y}, S^2 \times S^2, J_0 \right) \,=\, \min \left\{ 1-|x|, 1-|y| \right\} .
$$
(The corresponding $J_0$-holomorphic disc on $S^2 \times S^2$ with boundary on
$T_{x,y}$ can be seen in both $S^2 \times S^2$ and $\Box$, see the figure below.)

\begin{figure}[ht]
 \begin{center}
  \psfrag{mu}{$\mu$}
  \psfrag{t}{$\times$}
  \leavevmode\epsfbox{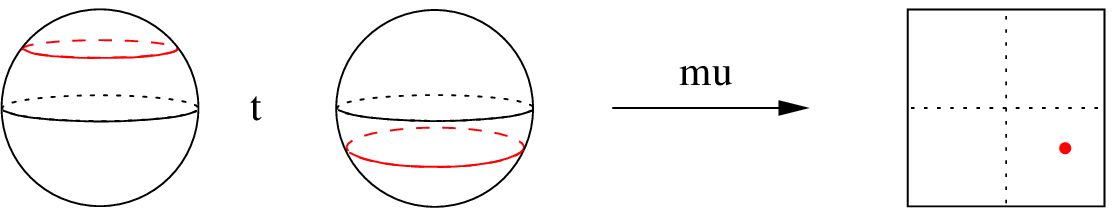}
 \end{center}
 \label{fig.mu}
\end{figure}
%
%

\ni
Therefore,
\begin{equation}  \label{e:xy}
S_{T^2}^e (x,y) \,=\, \min \left\{ 1-|x|, 1-|y| \right\}
                \,=\, 1+ \min \left\{ \pm \.x, \pm \.y \right\}
\end{equation}
if $(x,y) \neq (0,0)$.
In other words, on a punctured neighbourhood of $0 \in H_1(T^2;\RR)$ the function
$S_{T^2}^e -1$ is the minimum of {\bf four} linearly independent functionals.
We shall now show that on an open and dense set near $0 \in H_1(\Theta;\RR)$ the function
$S_{\Theta}^e -1$ is the minimum of {\bf three} linearly independent functionals.
These two functions are thus not related by a linear isomorphism of $\RR^2$.
Hence the tori $T^2$ and $\Theta$ are not symplectomorphic.

\m
\ni
{\bf Computation of $S_{\Theta}^e$.}
Recall that
$$
\Theta \,=\,
\left\{ \frac 1{\sqrt 2} \Bigl( \e^{2\pi\i \ga} \gg(t), \e^{-2\pi\i \ga} \gg(t) \Bigr)  \right\}
$$
where $\gg$ is a curve as in the figure.

\begin{figure}[ht]
 \begin{center}
  \leavevmode\epsfbox{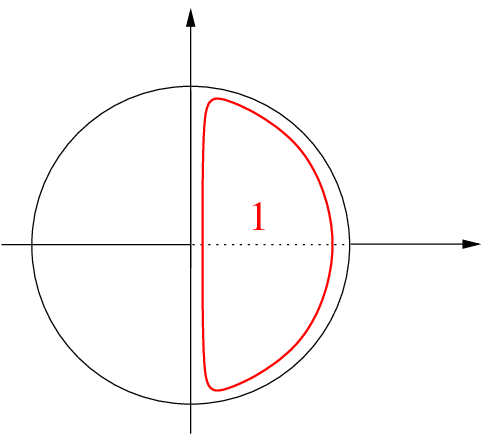}
 \end{center}
 \label{fig.gamma2}
\end{figure}
%
%

\ni
Note that $\Theta$ is invariant under the $S^1$-action
\begin{equation} \label{e:s1action}
(u,v) \mapsto \Bigl( \e^{2\pi\i \ga} u, \e^{-2\pi\i \ga} v \Bigr)
\end{equation}
and that
$$
\mu (\Theta) \,=\, \left\{ \pi \bigl( |\gg(t)|^2, |\gg(t)|^2 \bigr) \right\} \.=:\, \gs
$$
is a segment in the diagonal $\ell$ of $\Box$.

\begin{figure}[ht]
 \begin{center}
  \psfrag{1}{$1$}
  \psfrag{2}{$2$}
  \leavevmode\epsfbox{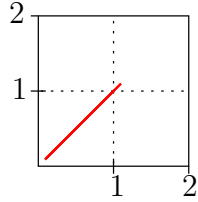}
 \end{center}
 \caption{The image $\mu(\Theta) \in \Box$ of the twist torus.}
 \label{fig.twist}
\end{figure}
%
%

As a basis of $H_1(\Theta;\ZZ)$ we take the class $[\Gamma]$
represented by the curve~$\Gamma$ in the diagonal~$\Delta$
and the class $[\tau]$ represented by an orbit of~\eqref{e:s1action}.
By the equivariant Weinstein Neighbourhood Theorem,
we can choose nearby Lagrangian tori $\Theta_{s,t}$ that are also invariant
under the action~\eqref{e:s1action}.
This means that for each such torus, $\mu (\Theta_{s,t})$ is a segment
$\gs_{s,t}$
parallel to~$\gs$.
The meaning of the deformation parameters $s$ and~$t$ is as follows.
Let $\gl$ be a primitive of the standard symplectic form~$\go$ on $\RR^4$ ($\go=d\lambda$).
Pick a curve~$\Gamma_{s,t}$
on the torus $\Theta_{s,t}$ that is close to the curve~$\Gamma$ on the torus $\Theta$.
Then $1+s$ is the integral of $\lambda$ over~$\Gamma_{s,t}$.
The parameter~$t$ is determined by the condition that
each $S^1$-orbit $\frak o$ on $\Theta_{s,t}$ has action $\int_{\frak o} \gl = t$.

\begin{figure}[ht]
 \begin{center}
  \psfrag{2}{$2$}
  \psfrag{s}{$\gs$}
  \psfrag{ss}{$\gs_{s,t}$}
  \psfrag{l}{$\ell$}
  \psfrag{ls}{$\ell_{s,t}$}
  \psfrag{t}{$t$}
  \psfrag{1x}{$1+x$}
  \psfrag{1y}{$1+y$}
  \leavevmode\epsfbox{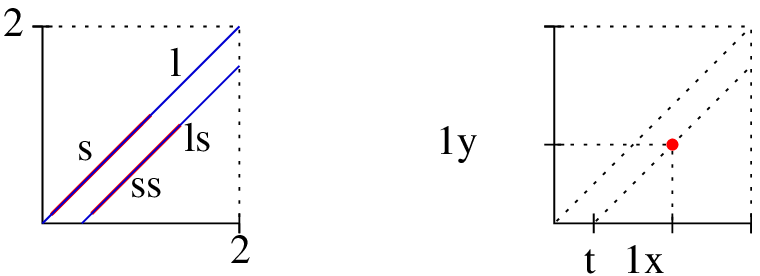}
 \end{center}
\end{figure}
%
%

In order to compute $e(\Theta_{s,t})$, we show that for $t \neq 0$ the torus
$\Theta_{s,t}$ is Hamiltonian isotopic to a product torus.
Fix $(s,t)$, and let $\ell_{s,t}$ be the half-open segment constructed
by intersecting the line containing $\gs_{s,t}$ with~$\Box$.
Let $x,y$ be such that $(1+x,1+y) \in \ell_{s,t}$.
Then both $T_{x,y}$ and $\Theta_{s,t}$ belong to the set
$\Sigma_{s,t} := \mu^{-1} \left(\ell_{s,t}\right)$.
Note that $T_{x,y} = T(1+x,1+y)$ is $S^1$-invariant,
and each $S^1$-orbit on this torus has action $(1+x)-(1+y) = x-y$.
On the other hand, all $S^1$-orbits contained in $\mu^{-1}(\ell_{s,t})$
have the same action. Therefore,
\begin{equation} \label{e:1}
t \,=\, x-y .
\end{equation}
Assume now that $t \neq 0$.
Then $\ell_{s,t}$ does not intersect the corner $(0,0)$ of $\Box$.
The $S^1$-action~\eqref{e:s1action} on $\Sigma_{s,t}$ is therefore free,
and $\Sigma_{s,t}$ smoothly splits as
\begin{equation}  \label{e:split}
\Sigma_{s,t} \,=\, D_{s,t} \times S^1
\end{equation}
where $D_{s,t} = \Sigma_{s,t} / S^1$ is a disc.
Denote by $\pi$  the projection $\Sigma_{s,t} \to D_{s,t}$.
There is an area form $\go_{s,t}$  on $D_{s,t}$ such that
$\go |_{\Sigma_{s,t}} = \pi^* \:\!\go_{s,t}$.
For each disc $D\subset \Sigma_{s,t}$
we thus have
\begin{equation} \label{e:symp}
\int_D \go \.=\, \int_{\pi(D)} \go_{s,t} .
\end{equation}
The sets $c_T := T_{x,y} / S^1$ and $c_\Theta := \Theta_{s,t} / S^1$
are smoothly embedded circles in~$D_{s,t}$.
The curve $\Gamma_{s,t}$ bounds a disc of symplectic area $1+s$;
moreover, such a disc can be found inside~$\Sigma_{s,t}$.
By~\eqref{e:symp}, the circle $c_\Theta=\pi (\Gamma_{s,t})$
also encloses symplectic area~$1+s$.
%
Again by~\eqref{e:symp},
the symplectic area enclosed by~$c_{T}$ equals the integral of $\gl$
over a lift of~$c_T$ to a circle in $T_{x,y}$ that is contractible in~$\Sigma_{s,t}$
(such a lift is unique up to homotopy).
One easily checks that the smaller of the two coordinate circles
whose product is $T_{x,y}$ is the required  lift.
Therefore, $c_T$ encloses symplectic area $1+\min \{x,y\}$.

Choose $(x,y)$ such that, in addition to~\eqref{e:1},
\begin{equation} \label{e:2}
s \.=\. \min \left\{ x,y \right\} .
\end{equation}
The two circles $c_T$ and $c_\Theta$ then enclose the same  area,
and hence are  Hamiltonian isotopic in~$D_{s,t}$.
The equivariant Hamiltonian lift of this isotopy to
$M_{\aff}$ yields a Hamiltonian isotopy from
$T(1+x,1+y) = T_{x,y}$ to $\Theta_{s,t}$.

Equations~\eqref{e:1} and \eqref{e:2} are equivalent to the set equality
$$
\left\{ x,y \right\} \,=\, \left\{ s, s+|t| \right\} .
$$
Since the displacement energy is invariant under symplectomorphisms,
$$
S^e_\Theta (s,t) \,=\, e \left( \Theta_{s,t} \right) \,=\,
e \left( T_{x,y} \right) \,=\, S_{T^2}^e (x,y)
$$
and hence, in view of~\eqref{e:xy}, for $t\ne0$ we have
\begin{eqnarray*}
S^e_\Theta (s,t) - 1 &=& \min \left\{ \pm s, \pm (s+|t|) \right\} \\
  &=& \min \left\{  s, -s-|t| \right\} \\
  &=& \min \left\{  s, -s \pm t \right\}  ,
\end{eqnarray*}
that is, away from a line the function $S^e_\Theta (s,t)-1$
is the minimum of three linearly independent functionals near $0$, as claimed.
\proofend

The idea of the above proof is summarized in the following figure.

\begin{figure}[h]
  \psfrag{S1}{$\Sigma_{s,0}/S^1$}
  \psfrag{S2}{$\Sigma_{s,t}/S^1$}
  \psfrag{T1}{$\Theta_{s,0}/S^1 \not\approx \T_{x,x}/S^1$}
  \psfrag{T2}{$\Theta_{s,t}/S^1 \approx \T_{x,y}/S^1$}
  \psfrag{l}{$\leadsto$}
    \includegraphics[width=9.0cm]{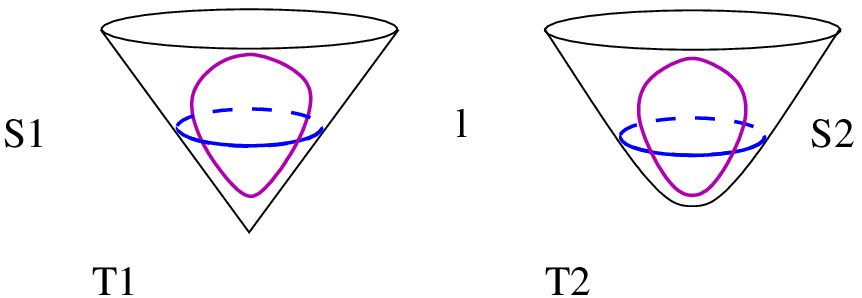}
\end{figure}

Assume now that $t=0$.
It has recently been shown in~\cite{FOOO3} that the non-monotone tori
$\Theta_{s,0}$ with $s>0$
near $\Theta$ are non-displaceable,
and hence in particular not product tori.
This can also be proved using pearl homology with Novikov coefficients
(cf.~the next section).
Assume now that $t=0$ and $s<0$, that is, the size of the disc in $\Delta$
cut out by $\Theta_{s,0}$ is less than~$1$.
Thus $\mu \bigl(\Theta_{s,0} \bigr)$ lies in the open lower left triangle of
$\Box$.
A rotation of the first factor of $S^2 \times S^2$ displaces $\Theta_{s,0}$,
whence $\Theta_{s,0}$ is displaceable.
Moreover, its displacement energy germ is $S^e_{\Theta_{s,0}}(s',t') =
1+s+s'-t'$,
whence $\Theta_{s,0}$ is not a product torus.
The set of neighbours of $\Theta = \Theta_{0,0}$ thus has the following structure.
\begin{eqnarray*}
\begin{array}{lcl}
\Theta_{0,0} &:& \phantom{non\;}\mbox{monotone twist torus, non-displaceable}, \\
\Theta_{s>0,0} &:& \mbox{non-monotone twist torus, non-displaceable}, \\
\Theta_{s<0,0} &:& \mbox{non-monotone twist torus, displaceable}, \\
\Theta_{s,t \neq 0} &:& \mbox{non-monotone product torus, displaceable}.
\end{array}
\end{eqnarray*}
\begin{figure}[ht]
 \begin{center}
  \psfrag{s}{$s$}
  \psfrag{t}{$t$}
  \leavevmode\epsfbox{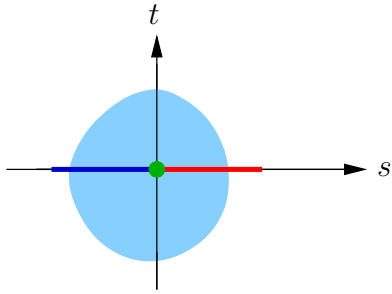}
 \end{center}
 \caption{Tori nearby $\Theta$.}
\end{figure}
%
%

\ni
A similar picture holds near each twist torus
in a product of complex projective spaces.

\section{Back to $\RR^{2n}$: inversion trick}

Recall that the many twist tori in $\RR^{2n}$ form only $n$ symplectic equivalence classes,
represented by
$$
T^n, \quad \Theta^1, \quad \Theta^2, \; \dots, \;\Theta^{n-1}
$$
or also by
$$
T^n, \quad \Theta \times T^{n-2}, \quad \Theta_1 (\Theta) \times T^{n-3},\; \dots, \;
\Theta_1 ( \Theta_1 ( \dots ))
$$
(see~\cite{Ch-tori}).
The reason why one does not get new tori by twisting more than once
is that in $\RR^{2n}$ there is enough room to ``untwist under the twist operation''.
This room is lacking in closed manifolds such as products of complex projective spaces,
and this is why our construction gives more different tori there.

However, one can construct many more different exotic tori in $\RR^{2n}$ by
performing the following inversion trick.
Notice that, by construction, each twist torus $T$ in $D^{2n}(1+\eps)\subset \RR^{2n}=\CC^n$
does not intersect the coordinate hyperplanes $\{z_m=0\}$.
Consider  $T$ as a monotone torus in
$$
 M = \CC P^{\ell_1}\times \dots\times \CC P^{\ell_j}.
$$
For each $k\in\{1,\dots,j\}$ remove from $\CC P^{\ell_k}$ one of the~$\ell_k$
coordinate hyperplanes.
The resulting manifold is a product of $j$ open balls.
This product of balls can be again put,
in the standard way, into~$\RR^{2n}$.
We thus construct a new monotone {\bf inverted torus}, $T'$ in~$\RR^{2n}$.
It turns out that inverted tori can be symplectically
different from twist tori.

For example, consider the torus $\Theta_1^1 \Theta^1 $ in~$\CC P^3$.
There are three ways to construct an inverted torus $T'$ out of it,
by removing one of the three coordinate hyperplanes.
By computing the displacement energy germs~$S^e_{T'}$,
one easily shows that
two of these three tori are symplectomorphic
neither to twist tori nor to each other.
These two tori are the only new monotone tori in $\RR^6$ that
can be constructed by inversion.
However, in higher dimensions there are many more possibilities
for performing inversion, and
the number of new tori in  $\RR^{2n}$ grows exponentially with~$n$.
The classification of inverted tori in lower dimensions
is still work in progress~\cite{CS1};
the complete classification in all dimensions does not look
feasible at the moment.


\section{Non-displaceability}

A subset $A$ of a symplectic manifold $(M,\go)$ is {\bf displaceable}
if there exists a Hamiltonian diffeomorphism $\Phi$ of $(M,\go)$ such that
$\Phi (A) \cap A = \varnothing$;
in other words, the displacement energy $e(A)$ of $A$ in $M$ is finite.

\m \ni
{\bf Theorem 2.}
{\it Twist tori in $\times_n S^2$ and $\CC P^n$ are not displaceable.}

\m
For the proof we use pearl (co)homology, as developed by Biran and Cornea~\cite{BC1,BC2}.
It is conceivable that one can also use Floer (co)homology as developed by
Fukaya, Oh, Ohta and Ono~\cite{FOOO1,FOOO2}.

Let $T$ be a twist torus.
Take a perfect Morse function $f \colon T \to \RR$, and let
$\Lambda := \ZZ [H_2(M,T)]$
be the group ring of the abelian group $H_2(M,T)$.
The pearl cochain complex $(C^*(f) \otimes \Lambda, d_\mathrm{P})$,
where $C^*(f)$ is the free abelian group
generated by the critical points of~$f$
and  $d_\mathrm{P}$ is a $\Lambda$-linear  differential of degree $+1$,
computes the pearl cohomology $\HP^* (T)$
(strictly speaking, only the $\ZZ_2$-version is written up at the moment,
but since tori admit a spin structure,  an appropriate coherent
orientation system should provide the signs required for the $\ZZ$-version).

For elements of $\Lambda$,
we use multiplicative notation in~$H_2(M,T)$,
writing $\widehat a$ for the multiplicative representation of
an element $a\in H_2(M,T)$.
The degree of $x\otimes \widehat D $ (which we abbreviate to $\widehat D\. x$)
is $|x| + \mu(D)$,
where $|x|$ is the Morse index of $x$ and $\mu(D)$ is the Maslov index of~$D\in H^2(M,T)$.

Take a generic Riemannian metric $g$ on $T$ and a generic almost complex structure~$J$ on~$M$.
The differential  $d_\mathrm{P}$
 is defined by counting not only gradient lines of~$f$ with respect to~$g$,
but also pearly trajectories formed by Morse lines and $J$-holomorphic discs
with boundary on~$T$ that pass through a given point in~$T$.
Then~$d_\mathrm{P}$ is of the form
$$
d_\mathrm{P} \,=\, d_0 + d_2 + d_4 + \dots
$$
where
$$
d_k\colon C^*(f)\otimes\Lambda \to C^{*-k+1}(f)\otimes\Lambda
$$
accounts for the pearly trajectories whose discs have Maslov index sum~$k$.
(See the figure below for an example of a pearly trajectory with one disc $\gs$ of Maslov index~$2$
that contributes to $d_2$: $\,d_2 (y)=\widehat{[\gs]}\, x   + \dots$.)

\begin{figure}[ht]
 \begin{center}
 \psfrag{x}{$x$}
 \psfrag{y}{$y$}
 \psfrag{s}{$\gs$}
 \psfrag{T}{$\Theta$}
  \leavevmode\epsfbox{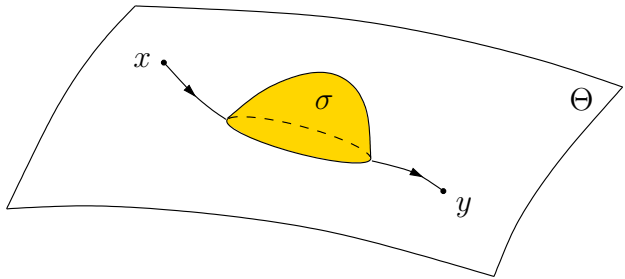}
 \end{center}
 \caption{A pearly trajectory.}
 \label{fig.pearl}
\end{figure}
%

Recall  that the pearl cohomology $\HP^*(T)$ is isomorphic to
the Floer cohomology of~$T$ with coefficients in~$\Lambda$,
and that $T$ is not displaceable in $M$ if the Floer cohomology does not
vanish~\cite{BC2}.
Therefore, in order to prove that $T$ is non-displaceable,
it suffices to show that $\HP^* (T) \neq 0$.

Since~$f$ is perfect, $d_0$ vanishes.
Therefore, we can identify $C^*(f)$ with  $H^*(T)$.
If the dimension of~$T$ is~$2$, then $d_k$ also vanish for $k \ge 4$, and we have
$$
d_\mathrm{P} \,=\, d_2 \colon H^*(T)\otimes\Lambda \to H^{*-1}(T)\otimes\Lambda.
$$
%
%
Let $\gs_1, \dots, \gs_\ell$ be the $J$-holomorphic discs of Maslov index~2
with boundary on $T$ that pass (in a non-degenerate way) through a given point in~$T$.

\begin{lemma}
The component $d_{\gs_k}$ of the differential
$d_2 = \sum_{k=1}^\ell d_{\gs_k}$
acts on $\alpha\in H^*(T)$ by
\begin{equation} \label{e:ds}
d_{\gs_k} \ga \,=\, \pm\.\widehat{\gs}_k \,\iiota_{\partial\gs_k} \ga
\end{equation}
where the sign is determined by an appropriate system of coherent orientations,
$\partial\gs_k \in H_1(T)$ denotes the homology class of the  boundary of
the disc $\gs_k$, and $\iiota$ is the contraction $H_1(T)\otimes H^*(T)  \to H^{*-1}(T)$.
\end{lemma}
The reason why this formula holds is that intersection in homology
translates to contraction in cohomology.

\m
In view of this lemma, in order to compute the pearl cohomology
of $T$, it is important to find, for some regular almost complex structure~$J$,
all $J$-holomorphic discs of Maslov index~$2$ with boundary on~$T$.
Enumerating such discs also allows to prove~Theorem~1.
For twist tori whose primitive components are product tori or twist tori of the form $\Theta^k$,
we compute all such discs with respect to the standard complex structure~$J_0$.
For twist tori involving iterated twists, we use different almost complex structures.

\m \ni
{\bf Computation of $J_0$-holomorphic discs with $\mu =2$ for $T = \Theta$.}
\ni
Consider the twist torus  $T=\Theta$ in $M=S^2\times S^2$.
For this example, we explain how to compute the $J_0$-discs of Maslov index~$2$
and how they can be used to show that the pearl cohomology of $\Theta$ does not vanish.
To avoid possible problems with determining the correct signs,
we use $\ZZ_2$-coefficients in all computations, with $\Lambda$
denoting now $ \ZZ_2 [H_2(M,T)]$.

To find $J_0$-holomorphic discs with $\mu =2$, we first compute the candidate classes
in $H_2(M,\Theta)$ that may contain such a disc.
Denote by $M_{\operatorname{aff}}$ the affine part
$(S^2 \setminus \infty) \times (S^2 \setminus \infty)$ of~$M$.
We choose
a basis $[D_{\lou{\Gamma}}], [D_\tau], [S_1], [S_2]$ of $H_2(M,\Theta)$ as follows.
First,
$D_{\lou{\Gamma}}$ is the disc in the diagonal $\Delta \subset M_{\aff}$ with boundary
$\Gamma = \bigl\{ \frac{1}{\sqrt 2} \bigl( \gg(t), \gg(t) \bigr) \bigr\}$.
Note that $\Gamma$ is one of the two components of the intersection of $\Theta$ with~$\Delta$.
Second, pick $(v,v) \in \Gamma$.
Then $D_\tau$ is  a disc in $M_{\aff}$ with boundary the orbit
$\tau = \left\{ \bigl( \e^{\i \ga} v, \e^{-\i \ga} v \bigr) \right\}$.
Finally, let $S_1=S^2 \times p$ and $S_2=p \times S^2$
for some point $p \in S^2\setminus\{0,\infty\}$.

For a closed holomorphic curve $\Sigma$ in $M$ disjoint from~$\Theta$,
and an oriented surface $D\subset M$ with boundary on~$T$,
the intersection number $\Sigma \cdot D$  is well defined.
Assume now that $D_0$ is a
$J_0$-holomorphic disc of Maslov index~$2$ with boundary on~$\Theta$.
Then $\Sigma \cdot D_0 \ge 0$ by positivity of intersection.
We now take for $\Sigma$ one of the five holomorphic curves
in the leftmost column of the table below.
Here, the complex number $w$ defining the degree~2 curve
$\left\{ z_1 z_2=w^2 \right\}$
is chosen in such a way that
the point $(w,w)$ lies in the interior of the disc~$D_{\lou{\Gamma}}$.
These curves are indeed disjoint from~$\Theta$.
The table presents the intersection number of these curves
with the cycles $D_{\Gamma}$, $D_{\tau}$, $S_1$, $S_2$;
the lowest row of the table gives the values of the Maslov class for these cycles.
\b
\begin{center}
  \begin{tabular}{  l | c  c    c  c  | l r}
                      & $D_\Gamma$  & $D_\tau$ & $S_1$ &  $S_2$ & $\phantom{-;}D_0$    \\
    \hline
    $S^2 \times 0$       & $0$  &  $-1$ & $0$  & $1$ & $b_2 \ge a_\tau$ &\quad (1) \\
    $S^2 \times \infty$  & $0$  &  $0$  & $0$  & $1$ & $b_2 \ge 0$ &\quad (2)  \\
    $0 \times S^2$       & $0$  &  $1$  & $1$  & $0$ & $b_1 \ge -a_\tau$ &\quad (3) \\
    $\infty \times S^2$  & $0$  &  $0$  & $1$  & $0$ & $b_1 \ge 0$ &\quad (4)  \\
    $\?\?\{z_1 z_2 = w^2\}$  & $1$  &  $0$  & $1$  & $1$ & $a_\Gamma+b_1+b_2 \ge 0$ &\quad (5)  \\
    $\.\mu$                & $2$  &  $0$  & $4$  & $4$ & $1=a_\Gamma+2(b_1+b_2)$ &\quad (6)
\end{tabular}
\end{center}

\b
\ni
Let
 $[D_0] = a_{\Gamma} [D_{\lou{\Gamma}}] + a_\tau [D_\tau] + b_1 S_1 + b_2 S_2$.

Positivity of intersections yields the inequalities (1)--(5)
for the coefficients $a_{\lou{\Gamma}}, a_\tau, b_1, b_2$.
The condition $\mu(D_0) = 2$ gives $2a_{\lou{\Gamma}}+4b_1+4b_2=2$, i.e.,
$1=a_{\lou{\Gamma}}+2(b_1+b_2)$.

\m
\ni
Subtracting (5) from (6) gives $1 \ge b_1+b_2 \ge 0$; thus, by (2) and (4), we must have
$$
(b_1,b_2) \,\in\, \left\{ (0,0), (1,0), (0,1) \right\} .
$$
If $(b_1,b_2)=(0,0)$, then (1) and (3) give $a_\tau=0$ and (6) gives $a_{\lou{\Gamma}}=1$,
that is $[D_0]=[D_{\lou{\Gamma}}]$.
In the same way we find the four other candidate classes
\begin{eqnarray*}
b_1=1, \quad b_2=0, \quad a_{\lou{\Gamma}}= -1, \quad a_\tau \in \{-1,0\} , \\
b_1=0, \quad b_2=1, \quad a_{\lou{\Gamma}}= -1, \quad a_\tau \in \{0,1\} .\phantom{-\,}
\end{eqnarray*}
We therefore have the five candidate classes
\begin{eqnarray*}
\phantom{-} \.[D_{\lou{\Gamma}}] &&\\
-\.[D_{\lou{\Gamma}}] &\?- \.[D_\tau] &+  \.[S_1] \\
-\.[D_{\lou{\Gamma}}] &\phantom{- \.[D_\tau]} & \.+  \.[S_1] \\
-\.[D_{\lou{\Gamma}}] &\phantom{- \.[D_\tau]} &\qquad\, +  \.[S_2] \\
-\.[D_{\lou{\Gamma}}] &+ \.[D_\tau] & \qquad\, +  \.[S_2]
\end{eqnarray*}

\begin{lemma} \label{l:five}
For each point $u\in\Theta$, each of the above five classes in $H_2(M,\Theta)$
is represented by  a unique
$J_0$-holomorphic disc $D_0$ such that $u\in \partial\. D_0\subset \Theta$.
\end{lemma}

The proof uses only complex analysis in dimension~1,
cf.~\cite{Au}.

\begin{lemma}
The $J_0$-holomorphic discs of Lemma~\ref{l:five}  are regular.
\end{lemma}

For the proof we use the holomorphic $S^1$-action.

\b \ni
{\bf Non-vanishing of $\.\HP^*(\Theta)$.}
\ni
We write $C^\ell (\Theta) $ for $H^\ell(\Theta) \otimes \Lambda$.
With respect to this grading, the differential $d_2$ has degree~$-\.1$.
We show that $\HP^0 (\Theta) \neq 0$, that is,
$ d_2 (H^1(\Theta) \otimes \Lambda)$
is not all of $H^0(\Theta) \otimes \Lambda$.
Since $H^0(\Theta)$ is canonically isomorphic to~$\ZZ$,
we can identify
$H^0(\Theta) \otimes \Lambda$ with~$\Lambda$.
It thus suffices to show that
$$
1 \notin d_2 \left( H^1 (\Theta) \otimes \Lambda \right) .
$$
Abusing notation, write  for the elements of $H_1(\Theta)$
represented by $\Gamma$ and~$\tau$.
This elements form a basis of~$H_1(\Theta)$.
Let $\Gamma^*, \tau^*$ be the dual basis of~$H^1(\Theta)$.
Consider the generators of $H_2(M,\Theta)$ given, in multiplicative notation,
by $R=\widehat{[D_\lou{\Gamma}]},$ $T=\widehat{[D_{\tau}]}$,
and also (abusing notation again) the generators
$S_1=\widehat{S_1},$ $S_2=\widehat{S_2}$.
According to~\eqref{e:ds}, we have
\begin{eqnarray*}
\begin{array}{lllll}
d_{\gs_1} &=& R \,\iiota_{\.\Gamma} \\
d_{\gs_2} &=& R^{-1} T^{-1} S_1 \,\iiota_{-\Gamma-\.\tau} \\
d_{\gs_3} &=& R^{-1} S_1 \,\iiota_{-\Gamma} \\
d_{\gs_4} &=& R^{-1} S_2 \,\iiota_{-\Gamma} \\
d_{\gs_5} &=& R^{-1} T S_2 \,\iiota_{-\Gamma+\.\tau} .
\end{array}
\end{eqnarray*}
Since we work over $\ZZ_2$, this implies
\begin{eqnarray*}
\begin{array}{lllll}
d_2 \Gamma^* &=& \sum_{k=1}^5 d_{\gs_k} \Gamma^* &=& R+R^{-1} \left( T^{-1} S_1 + S_1 + S_2 + T S_2  \right) \\
d_2 \tau^* &=& \sum_{k=1}^5 d_{\gs_k} \tau^* &=&  R^{-1} \left( T^{-1} S_1 + T S_2 \right) .
\end{array}
\end{eqnarray*}
Let $\gf \colon \Lambda \to \ZZ_2 [R,R^{-1}]$ be the  epimorphism defined by
$$
\gf (R) = \gf (T) = \gf (S_2) = R, \quad \gf (S_1) = 1 .
$$
Then
\begin{eqnarray*}
\begin{array}{lll}
\gf (d_2 \Gamma^*) &=& R+R^{-1} \left( R^{-1} +1+R+R^2 \right) \;=\; R^{-2} \left( 1 + R + R^2 \right) ,\\
\gf (d_2 \tau^*) &=& R^{-1} \left( R^{-1} + R^2 \right) \;=\; R^{-2} \left( R^3+1 \right) \;=\;
R^{-2} (R+1) (R^2+R+1).
\end{array}
\end{eqnarray*}
Thus $\gf$ maps $d_2 \bigl( H^1 (\Theta) \otimes \Lambda \bigr)$ onto the proper ideal $(R^2+R+1)$
in $\ZZ_2[R,R^{-1}]$.
On the other hand, if $1 \in d_2 \bigl( H^1(\Theta) \otimes \Lambda \bigr)$,
then $\gf \bigl( d_2 \bigl( H^1(\Theta) \otimes \Lambda \bigr) \bigr) = \ZZ_2[R,R^{-1}]$.
This contradiction proves that $\Theta$ is non-displaceable.

\begin{remarks}
{\rm
{\bf \;\;1.}
It is important for the proof that we distinguish the elements in $H_2(M,\Theta)$,
that is, work with the full ring $\Lambda = \ZZ_2 [H_2(M,\Theta)]$.
If instead we work over the smaller ring of Laurent polynomials $\ZZ_2[t,t^{-1}]$,
obtained from $\Lambda$ by mapping $A$ to $t^{\mu(A)/2}$
(that is, $R \mapsto t$, $T \mapsto 1$, $S_1 \mapsto t^2$, $S_2 \mapsto t^2$),
then $d_2 \Gamma^* = t+t^{-1}(t^2+t^2+t^2+t^2) = t$,
and so $d_2 (t^{-1} \otimes \Gamma^*) =1$, whence $\HP^0(\Theta)=0$
over this coefficient ring.

\m \ni
{\bf 2.}
Proceeding in the same way, we prove that
$\HP^0 \left( \Theta^{n-1}, \times_n S^2 ;\ZZ_2\right) \neq 0$
for all $n$ even.
For $n$ odd, one has $\HP^0 \left( \Theta^{n-1}, \times_n S^2 ;\ZZ_2\right) = 0$,
but $\HP^0 \left( \Theta^{n-1}, \times_n S^2 ;\ZZ\right)$ is non-trivial.

The same holds for twist tori $\Theta^{n-1}$ in complex projective space $\CP^n$,
where working over $\ZZ_2$ suffices for $n$ odd, while one needs to work over $\ZZ$ for $n$ even.
}
\end{remarks}

\m \ni
{\bf Non-displaceability in higher dimensions.}
\ni
To prove that a general twist torus $T$ is non-displaceable,
we also show that the zero degree component
of the pearl cohomology $\HP^* ( T)$ does not vanish.
The task of proving this seems much more complicated because,
in general, the operators $d_4, d_6, \dots$ contributing to
the differential $d_\mathrm{P}$ do not have to vanish and are hard to compute.
It turns out, however, that it suffices to consider
only the term~$d_2$.
As in the case $T=\Theta$ we use the grading on $H^*(T)\otimes\Lambda$
with respect to which $\deg \bigl( H^*(T)\otimes\Lambda \bigr) = *$.
The operator $d_k$ then has degree~$1-k$.
This implies $d_2^2=0$.
The following lemma can be seen as an immediate application of the
spectral sequence theory:

\begin{lemma}
If $H^0 \bigl( H^*(T)\otimes\Lambda, d_2 \bigr) \ne 0$ and
$H^\ell \bigl( H^*(T) \otimes\Lambda, d_2 \bigr) = 0$ for $\ell>0$,
then $\HP^* (T) \cong H^*\bigl( H^*(T) \otimes \Lambda, d_2 \bigr)$.
\end{lemma}

In view of this lemma, it suffices to show that the cohomology
of $d_2$ does not vanish in degree~$0$
and vanishes in all other degrees.

Pick free generators $S_1,\dots,S_j$ of~$H_2(M)$, and
pick $R_1,\dots,R_n \in H_2(M,T)$ such that $S_1,\dots,S_j, R_1,\dots,R_n$
are free generators of~$H_2(M,T)$ (written multiplicatively).
Then, in particular, $n$ is the dimension of~$T$.
Denote by $\delta$ the boundary map $H_2(M,T)\to H_1(T)$.
Let $q_1, \dots, q_n$ be the basis of~$H^1(T)$ (additively written)
dual to the basis $\delta R_1,\dots,\delta R_n$ of $H_1(T)$.
Consider the {\bf potential function}
$$
          U:=\, \sum_{k=1}^\ell \pm\.\.\widehat {[\gs_k]},
$$
where $\gs_1,\dots ,\gs_\ell$ are  the $J$-holomorphic discs with Maslov index $2$
passing generically through a point in~$T$.
Then~\eqref{e:ds} translates into
\begin{equation}
d_2 \;\!\ga \,=\, \sum_{k=1}^m \, R_k \,\frac{\partial \.U}{\partial R_{\.k}}\,
           \iiota_{q_{k}} \ga.
\end{equation}
Therefore, the complex $\bigl( H^*(T)\otimes \Lambda, d_2 \bigr)$
is the {\bf Koszul complex} \cite{Ei,Wei} 
associated with the toric differential
$$
\Bigr( R_1 \,\frac{\partial\.U}{\partial R_1\?\?},\dots, R_n \,\frac{\partial\.U}{\partial R_n\?\?}\.\Bigl)
\,\,=:(v_1,\dots,v_n)
$$
of the Laurent polynomial~$U$.

By Koszul's Theorem 
(\cite[Corollary 17.5]{Ei} or \cite[Corollary 4.5.4]{Wei}), 
in order to prove that
$H^\ell \bigl( H^*(T) \otimes \Lambda, d_2 \bigr)$ vanishes exactly for $\ell>0$,
it suffices to show that the elements $v_1,\dots,v_n$ (possibly after a reordering)
form a {\bf regular sequence}, that is, each $v_k$ is not a zero divisor
in the quotient ring $\Lambda/(v_1,\dots,v_{k-\?1})$.
A convenient way to prove the regularity property is by using
a homomorphism $\varphi$ from $\Lambda$ to another, less complicated, ring~$A$.
If the sequence $\varphi(v_1),\dots,\varphi(v_n)$ is regular in~$A$,
then $v_1,\dots,v_n$ is regular in~$\Lambda$.
We already used this homomorphism trick
(without mentioning regular sequences)
in the argument for~$T=\Theta$.

A nice example of such a homomorphism is as follows:
consider the homomorphism $\varphi$ from $\Lambda$
to the ring of Laurent polynomials in $n$ complex variables $z_1,\dots, z_n$,
such that each generator $S_\ell$ is mapped to a nonzero complex number~$c_\ell$,
and each generator $R_k$ is mapped to the monomial~$z_k$.
Then the sequence $\varphi(v_1),\dots,\varphi(v_n)$ is
regular if and only if all critical points of the function~$\varphi(U)$
are isolated.

\section{Other constructions of exotic tori in dimension four}

Different constructions of an exotic torus in $S^2 \times S^2$ were given
by Biran--Cornea~\cite{BC2},
Entov--Polterovich~\cite{EP},
Fukaya--Oh--Ohta--Ono~\cite{FOOO3}.
Also, there is an exotic torus in $S^2 \times S^2$
coming from the geodesic flow on $T^*S^2$, see Albers--Frauenfelder~\cite{AF}.
%
It is known that these tori are all Hamiltonian equivalent.
Agn\`es Gadbled~\cite{Ga} is on the way to show that the twist torus $\Theta$
is also equivalent to these tori.
She has already shown that~$\Theta$ in~$\CP^2$
is Hamiltonian equivalent to the torus in $\CP^2$ described by Biran--Cornea~\cite{BC2}.

\section{The twist torus in blow-ups of $\CP^2$}

Let $\CP^2$ be complex projective space endowed with the standard symplectic form~$\go$
normalized such that $\int_{\CP^1}\go = 3$.
Then $\Theta \subset B^4(3) = \CP^2 \setminus \CP^1$, and $\Theta$ is monotone in $\CP^2$.
Versal deformations show that $\Theta$ is not symplectomorphic to the Clifford torus.
The same holds true in the monotone blow-up of $\CP^2$ at one and two points, see the figure.
The segment $\mu (\Theta)$ is just too long to fit into the blow up of $\CP^2$ at three points.
Can our construction of $\Theta$ be modified so that one obtains an exotic monotone torus in
the blow-up of $\CP^2$ at three points?

\begin{figure}[h]
 \begin{center}
  \leavevmode\epsfbox{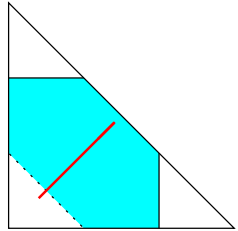}
 \end{center}
 \caption{Does $\Theta$ fit into $\CP^2$ blown-up at three points?}
 \label{fig.blow}
\end{figure}
%

\enddocument